\documentclass[onefignum,onetabnum]{siuro210301}
\usepackage{lipsum}
\usepackage{amsfonts}
\usepackage{graphicx}
\usepackage{epstopdf}
\usepackage{algorithmic}
\usepackage{booktabs}
\usepackage{multirow}
\usepackage{caption}
\usepackage{subcaption}
\usepackage{comment}

\Crefname{ALC@unique}{Line}{Lines} 

\newcommand{\R}{\mathbb{R}}

\DeclareMathOperator*{\argmin}{arg\,min}                   
\renewcommand{\t} {^{\top}}                                

\renewcommand{\phi}{\mathbf{\varphi}}

\newcommand{\bfd}{\mathbf{d}}

\newcommand{\bfF}{\mathbf{F}}

\newcommand{\bfA}{\mathbf{A}}

\newcommand{\bfE}{\mathbf{E}}
\newcommand{\bfU}{\mathbf{U}}

\newcommand{\bfR}{\mathbf{R}}
\newcommand{\bfI}{\mathbf{I}}

\newcommand{\bfb}{\mathbf{b}}

\newcommand{\bfB}{\mathbf{B}}

\newcommand{\bfx}{\mathbf{x}}

\newcommand{\bfe}{\mathbf{e}}
\newcommand{\bfu}{\mathbf{u}}
\newcommand{\bfy}{\mathbf{y}}

\newcommand{\bfM}{\mathbf{M}}

\newcommand{\bfr}{\mathbf{r}}
\newcommand{\bfs}{\mathbf{s}}
\newcommand{\bfv}{\mathbf{v}}
\newcommand{\bfV}{\mathbf{V}}

\newcommand{\bfQ}{\mathbf{Q}}

\newcommand\norm[1]{\left\lVert#1\right\rVert}

\ifpdf
  \DeclareGraphicsExtensions{.eps,.pdf,.png,.jpg}
\else
  \DeclareGraphicsExtensions{.eps}
\fi

\usepackage{enumitem}
\setlist[enumerate]{leftmargin=.5in}
\setlist[itemize]{leftmargin=.5in}

\newsiamremark{remark}{Remark}
\newsiamremark{hypothesis}{Hypothesis}
\crefname{hypothesis}{Hypothesis}{Hypotheses}
\newsiamthm{claim}{Claim}

\headers{Inexact genGK}{Bu}

\title{Inexact Generalized Golub-Kahan Methods for Large-Scale Bayesian Inverse Problems}

\author{Yutong Bu\thanks{Emory University \textbf{}
  (\email{ybu4@emory.edu})}}

\dedication{\small\textit{Project advisor: Julianne Chung\thanks{Department of Mathematics, Emory University (\email{jmchung@emory.edu}).}}}

\usepackage{amsopn}

\ifpdf
\hypersetup{
  pdftitle={An Example SIURO Article},
  pdfauthor={D. Doe, P. T. Frank, and J. E. Smith}
}
\fi

\begin{document}

\maketitle

\begin{abstract}
Solving large-scale Bayesian inverse problems presents significant challenges, particularly when the exact (discretized) forward operator is unavailable. These challenges often arise in image processing tasks due to unknown defects in the forward process that may result in varying degrees of inexactness in the forward model.  Moreover, for many large-scale problems, computing the square root or inverse of the prior covariance matrix is infeasible such as when the covariance kernel is defined on irregular grids or is accessible only through matrix-vector products. This paper introduces an efficient approach by developing an inexact generalized Golub-Kahan decomposition that can incorporate varying degrees of inexactness in the forward model to solve large-scale generalized Tikhonov regularized problems. Further, a hybrid iterative projection scheme is developed to automatically select Tikhonov regularization parameters. Numerical experiments on simulated tomography reconstructions demonstrate the stability and effectiveness of this novel hybrid approach.
\end{abstract}



\section{Introduction}
\label{sec:intro}

Inverse problems are prevalent in many scientific applications, but they are usually hard to solve due to their large scale. Inverse problems are typically ill-posed, whereby the number of unknown parameters can be significantly larger than the size of the observed datasets. This projects aims to provide an efficient method to solve large-scale inverse problem through a Bayesian approach, specifically when the exact forward model is not available, as described in \cite{simoncini2003theory, elfving2018unmatched}, but where an approximate model can be obtained. 

Consider a linear inverse problem of the form 
\begin{equation}
    \mathbf{As+\boldsymbol{\epsilon} = d}, 
    \label{eq:inverse problem}
\end{equation}
where $\mathbf{d} \in \mathbb{R}^m $ is the observed data, $\mathbf{A}\in\mathbb{R}^{m\times n} $ is an ill-posed matrix representing the discretized forward model, $\mathbf{s} \in \mathbb{R}^n$ is the desired parameters, and $\boldsymbol{\epsilon}\in\mathbb{R}^m$ is the additive Gaussian noise in the data. Assume $\mathbf{A}$ is only accessible through matrix-vector products (MVPs) with some potential inexactness involved, $\boldsymbol{\epsilon} \sim \mathcal{N}(\mathbf{0, R)}$ where $\mathbf{R}$ is a symmetric positive definite (SPD) matrix with inexpensive inverse and square root, for example a diagonal matrix with positive entries.  The goal is to compute an approximation for $\mathbf{s\text{ given } A \text{ and }b}$.

Given the ill-posed nature of the problem from both the ill-conditioned $\bfA$ and the additive noise $\boldsymbol{\epsilon}$ in observation $\bfd$, a small error in $\bfd$ could lead to large deviations in the computed solution from the desired solution $\bfs$. Thus, a necessary regularization technique is required to obtain a meaningful recovery or approximation to the true solution. The Bayesian framework provides an effective approach for this type of problems, as it naturally incorporates regularization through prior knowledge, meanwhile, quantifies the uncertainty in the desired parameters $\mathbf{s}$ utilizing the posterior density function \cite{waqar_tutorial_2023, ChungGeneralized, calvetti2023bayesian}.

Following the Bayesian approach, the solution of the inverse problem could be represented by the probability distribution of $\mathbf{s \text{ given } d}$, denoted as $\mathbf{s|d}$. Assume $\mathbf{s} \sim \mathcal{N}(\boldsymbol{\mu}, \lambda^{-2}\mathbf{Q})$ as the prior of which $\mathbf{s}$ is a Gaussian random variable with mean $\boldsymbol{\mu}$ and SPD  covariance matrix $\mathbf{Q}$ and $\lambda^2$ is a scaling parameter for the precision matrix. Therefore, by Bayes’ Theorem, the posterior probability distribution function is

\begin{equation}
    p(\mathbf{s|d})  \propto p(\mathbf{d|s})p(\mathbf{s}) = \exp \left(-\frac{1}{2}\norm{\mathbf{As-d}}^2_{\mathbf{R}^{-1}} -\frac{\lambda^2}{2}\norm{\mathbf{s}-\boldsymbol{\mu}}^2_{\mathbf{Q^{-1}}}\right)
\end{equation}
where $\mathbf{\norm{x}_M=\sqrt{x^\top M x} }$ is a vector norm for any SPD matrix $\mathbf{M}$. Similar to the maximum likelihood estimator, the maximum a posteriori (MAP) estimate provides a  solution to \cref{eq:inverse problem} and can be obtained by minimizing the negative logarithm of the posterior
\begin{equation}
    \begin{aligned}
    \mathbf{s}_\lambda &=   \underset{\mathbf{s}}{\arg\min} -\log p(\mathbf{s|d})  \\
    &= \underset{\mathbf{s}}{\arg\min} \frac{1}{2}\norm{\mathbf{As-d}}^2_{\mathbf{R}^{-1}} +\frac{\lambda^2}{2}\norm{\mathbf{s-\mu}}^2_{\mathbf{Q}^{-1}}.
    \end{aligned}\label{eq:MAP}
\end{equation} 
In fact, $\bfs_\lambda$ is the solution of a general-form Tikhonov problem, for which many approaches have been developed to compute a reliable solution through hybrid iterative methods  \cite{ChungGeneralized, reichel2012tikhonov, kilmer2007projection, gazzola2014generalized}. However, in this work, we further account for the uncertainties in the forward operator, by assuming that exact MVPs with $\bfA$ may not be available.  Although such assumptions were considered for inexact Krylov methods \cite{simoncini2003theory}, they have not been considered in previous methods for solving generalized Tikhonov problems.

This work addresses the scenario when the forward operator $\bfA$ is not fully known. This is a common issue in signal and image processing tasks \cite{gazzola2021regularization}, such as device calibration \cite{chung2010efficient, golub2003separable}, blind deconvolution \cite{ji2011robust}, and super-resolution \cite{chung2006numerical}. In these cases, the MVPs with $\mathbf{A}$ and $\mathbf{A}^\top$ cannot be performed exactly, i.e. they are only available as approximations denoted as $\widehat{\bfA}$ and $\widetilde{\bfA}^\top$.

Although the errors in $\bfA$ can be interpreted as model error, and there is prior work on statistical approaches to handle model error uncertainty \cite{smith2024uncertainty}, interpreting such model errors as random variables in a Bayesian framework can be computationally infeasible for our problems of interest.  Another approach is to parameterize the forward model, hence reducing the number of unknown parameters defining the forward model \cite{chung2010efficient}. However, computational methods that can incorporate uncertainty arising from the forward model focus mainly on computing point estimates, rather than uncertainty. Here, we assess the forward operator's inexactness through error matrices instead of as random variables or parameterized models. The primary reason for this approach is due to its large scale ($\bfA \in \mathbb{R}^{m\times n}$), which leads to the creation of too many unknown parameters and becomes computationally intractable. Moreover, our approach can handle scenarios where the prior covariance matrix is too large for obtaining a factorization.

\paragraph{Main contributions}
In this paper, we propose the inexact generalized Golub-Kahan hybrid method to compute a solution for \cref{eq:MAP} where $\bfA$ contains errors. In particular, after an appropriate change of variables, the inexact generalized Golub-Kahan bidiagonalization methods is used to address inexactness in the forward model, meanwhile maintaining general to a rich class of covariance kernels. Further, a hybrid approach is adopted to automatically select a regularization parameter in solving the projected problem at each iteration.  Numerical experiments show that the inexact generalized hybrid method can achieve a solution with high accuracy that is comparable to its exact counterpart.

The paper is organized as follows. In \Cref{sec: change of variable}, we reformulate the problem through the change of variables to make it computationally feasible with the choice of covariance kernel $\bfQ$. \Cref{sec:background} reviews the existing related iterative methods: Golub-Kahan Bidiagonalization for solving least squares (LS) problems, inexact Golub-Kahan decomposition for solving LS problems with inexact $\bfA$ and $\bfA^\top$, and generalized Golub-Kahan Bidiagonalization for solving generalized LS problems. \Cref{sec:Inexact genGK} introduces the proposed inexact generalized Golub-Kahan  method which is further adapted to the hybrid scheme in solving the Tikhonov regularized problem. \Cref{sec:numerical experiment} presents various numerical experiments, and \Cref{sec:conclusion} concludes this paper with some remarks and future directions.

\section{Change of Variables}
\label{sec: change of variable}
Our goal is to compute the MAP estimate \cref{eq:MAP} derived from the Bayesian approach. This section presents a way to reformulate the problem to be more computationally feasible through a change of variables.

Equivalently, the above MAP estimate can be written as the standard-form Tikhonov problem where iterative methods have been developed to solve 
\begin{equation}
\min_\mathbf{s} \frac{1}{2}\norm{\mathbf{L_R(As-d)}}^2_2 +\frac{\lambda^2}{2}\norm{\bf{L_Q}(\bfs-\boldsymbol{\mu})}^2_2
\end{equation}
where $\mathbf{R}^{-1} = \mathbf{L_R}^\top\mathbf{L_R}$ and $\mathbf{Q}^{-1} = \mathbf{L_Q}^\top\mathbf{L_Q}$ are any symmetric matrix factorization (e.g., Cholesky). However, with an interest on incorporating Gaussian kernels in the prior for this problem, the computation of $\mathbf{Q^{-1} ,\ \mathbf{L_Q}^\top,\text{ or }\mathbf{L_Q}}$ could be very expensive. Thus, to avoid such computations while directly solving the Tikhonov normal equation or through priorconditioning, a change of variables  
\[
\mathbf{x} \leftarrow \mathbf{Q}^{-1} (\mathbf{s}-\boldsymbol{\mu}),\quad \mathbf{b}\leftarrow \mathbf{d-A}\boldsymbol{\mu} 
\]
could be applied \cite{ChungGeneralized}.

Then, we obtain an equivalent problem
\begin{equation}
    \mathbf{x}_\lambda  = \argmin_\mathbf{x} \frac{1}{2}\norm{\mathbf{AQx-b}}^2_{\mathbf{R}^{-1}} +\frac{\lambda^2}{2}\norm{\mathbf{x}}^2_\mathbf{Q} 
    \label{eq:reformulated MAP}
\end{equation}
with the MAP estimate $\mathbf{s}_\lambda = \boldsymbol{\mu}+\mathbf{Qx}_\lambda$ and the new Bayesian interpretation 
\begin{equation*}
\bf{b|x} \sim \mathcal{N}(\bf{AQx, R}^{-1}),\quad \bf{x}\sim\mathcal{N}(\bf{0},\lambda^{-2}\bf{Q}^{-1}).
\end{equation*}

Thus, the MAP estimate could be reformulated as a LS problem as in \cref{eq:reformulated MAP}, with $\mathbf{x}=\mathbf{Q}^{-1} (\mathbf{s}-\boldsymbol{\mu}),\ \mathbf{b}=\mathbf{d-A}\boldsymbol{\mu} $, $\mathbf{A}$ being a known forward process, and $\mathbf{Q}$ is a covariance matrix from the Matérn family (detailed explanation for the choice of $\bfQ$ can be found in \cref{sec:cov kernel}). $\mathbf{A}$ and $\mathbf{Q}$ are really large so they cannot be explicitly stored and are only accessible through MVPs as function handles.

\section{Background on Iterative Methods for Inverse Problems}
\label{sec:background}
The goal of this section is to introduce several established iterative methods in solving related problems, upon which the proposed inexact gen-GK bidiagonalization detailed in \Cref{sec:Inexact genGK} is built upon.

\subsection{Golub-Kahan Bidiagonalization}
Given an unregularized standard LS problem,
\begin{equation}
\min_{\bfx\in\mathbb{R}^n} \norm{\mathbf{Ax-b}}^2_2,
    \label{eq:standard setup}
\end{equation}
the Golub-Kahan (GK) bidiagonalization is one of the most common iterative methods which project the stated problem onto subspaces of increasing dimension (e.g. Krylov subspaces) \cite{golub1965calculating, chung2024computational}. It generates two sets of orthogonal vectors to span the Krylov subspaces $\mathcal{K}_k(\mathbf{AA^\top, b)} $ and $\mathcal{K}_k(\mathbf{A^\top A, A^\top b)}$. 

Given $\mathbf{A}\in\mathbb{R}^{m\times n} $ and $\mathbf{b}\in\mathbb{R}^m $, the GK process initiates $\mathbf{u}^{\text{s}}_1=\mathbf{b}/ \beta_1^{\text{s}}$ where $\beta^{\text{s}}_1=\norm{\mathbf{b}}$ and $\mathbf{v}^{\text{s}}_1=\mathbf{A^\top u}^{\text{s}}_1/ \alpha_1^{\text{s}}$ where $\alpha_1^{\text{s}} = \norm{\mathbf{A^\top u}^{\text{s}}_1}$. The superscript `s' denotes the basis generated for standard GK bidiagonalization. Then, at the $k$-th iteration, it method generates vectors $\mathbf{u}^{\text{s}}_{k+1}$ and $\mathbf{v}^{\text{s}}_{k+1} $ and scalars $\alpha^{\text{s}}_{k+1} $ and $\beta^{\text{s}}_{k+1}$, as diagonal and sub-diagonal entries in $\mathbf{B}_k^{\text{s}}$, through 
\begin{equation}
    \mathbf{u}^{\text{s}}_{k+1} = (\mathbf{Av}^{\text{s}}_k-\alpha^{\text{s}}_k\mathbf{u}^{\text{s}}_k)/\beta^{\text{s}}_{k+1},\qquad \mathbf{v}^{\text{s}}_{k+1}=(\mathbf{A^\top u}^{\text{s}}_{k+1}-\beta^{\text{s}}_{k+1}\mathbf{v}^{\text{s}}_k)/\alpha^{\text{s}}_{k+1},
\end{equation}
where the values of $\alpha^{\text{s}}_{k+1} $ and $\beta^{\text{s}}_{k+1}$ are chosen to ensure  $\norm{\mathbf{v}^{\text{s}}_{k+1} } =1$ and $\norm{\mathbf{u}^{\text{s}}_{k+1}}=1 $ respectively.

After $k$ iterations of the GK process, we obtain 
\[\mathbf{B}^{\text{s}}_k = \begin{bmatrix}
    \alpha^{\text{s}}_1 &0&\dots&0\\
    \beta^{\text{s}}_2 & \alpha^{\text{s}}_2&\ddots&\vdots  \\
    0& \beta^{\text{s}}_3 &\ddots&0\\
    \vdots&\ddots &\ddots & \alpha^{\text{s}}_k\\
    0 &\dots &0 & \beta^{\text{s}}_{k+1}
 \end{bmatrix}\in \mathbb{R}^{(k+1)\times k}, \quad \mathbf{U}^{\text{s}}_{k+1} = [\mathbf{u}^{\text{s}},\dots,\mathbf{u}^{\text{s}}_{k+1} ],\quad \mathbf{V}^{\text{s}}_k=[\mathbf{v}^{\text{s}}_1,\dots,\mathbf{v}^{\text{s}}_k],
 \]
 with following relationships 
 \begin{equation}
 \begin{aligned}
\mathbf{U}_{k+1}^{\text{s}}\beta_1^{\text{s}}\mathbf{e}_1 &= \mathbf{b}, \\
\mathbf{AV}^{\text{s}}_k&=\mathbf{U}^{\text{s}}_{k+1}\mathbf{B}^{\text{s}}_k,\\  \mathbf{A^\top U}^{\text{s}}_{k+1}&=\mathbf{V}^{\text{s}}_k(\mathbf{B}^{\text{s}}_k)^\top+ \alpha^{\text{s}}_{k+1}\mathbf{v}^{\text{s}}_{k+1}\mathbf{e}_{k+1}^\top,
     \label{eq:GKB matrix fac}
\end{aligned}
 \end{equation}
where $\bfe_i$ is the $i$th column of the identity matrix of appropriate size.

LSQR is an iterative method for solving \cref{eq:standard setup} where at the $k$th iteration, we seek a solution $\bfx_k^{\text{s}} \in \text{span}\{\bfV_k^{\text{s}} \}$.  Define the residual at $k$-th iteration as
\[
\bfr_k \equiv \bfA\bfx_k^{\text{s}} - \bfb = \bfA\bfV_k^{\text{s}}\bfy_k^{\text{s}} - \bfb = \bfU_{k+1}^{\text{s}}(\bfB_k^{\text{s}}\bfy_k^{\text{s}} - \beta_1^{\text{s}} \bfe_1),
\]
where $\bfy_k^{\text{s}} \in \R^k$. Then, we could solve for $\bfy_k^{\text{s}}$ by solving the following projected LS problem
\begin{equation}
\bfy_k^{\text{s}}  = \min_{\bfy \in \R^k} \norm{\bfB_k^{\text{s}} \bfy - \beta_1^{\text{s}} \bfe_1}_2^2,
\end{equation}
and recover the solution as $\bfx_k^{\text{s}} = \bfV_k^{\text{s}} \bfy_k^{\text{s}}$.

\subsection{Inexact Golub-Kahan Decomposition}
Consider the same linear system as in \cref{eq:standard setup}, but assume that the MVPs with $\mathbf{A\text{ and } A^\top}$ are only approximately available, i.e. at the $k$-th iteration, we have approximate MVPs,
\begin{equation}
\begin{aligned}
\widehat{\bfA}_k \bfx &= (\bfA + \bfE_k)\bfx \\
\widetilde{\bfA}_k^\top \bfy &= (\bfA + \bfF_k)^\top \bfy.
\end{aligned}
\label{eq:inexact A def}
\end{equation}
where $\bfE_k$ and $\bfF_k$ are the errors occured during the MVPs with $\bfA$ and $\bfA^\top$ at the $k$-th iteration. In this scenario, the inexact Golub-Kahan (iGK) decomposition \cite{gazzola2021regularization} provides an efficient method for computing a solution subspace.

The iteration-wise computation and matrix factorization of iGK are different from GK in that they adapt to the inexactness in MVPs. Similar to the notation used for GK, we use the superscripts i to represent the vectors and matrices computed using iGK. Initializing  $\mathbf{u}_1^{\text{i}}  = \mathbf{b}/\beta^{\text{i}} $ where $\beta^{\text{i}}=\norm{\mathbf{b}} $, and $\mathbf{v}_1^{\text{i}}=(\mathbf{A+F}_1)^\top \mathbf{u}_1^{\text{i}}/[\mathbf{L}^{\text{i}}]_{1,1} $ where $[\mathbf{L}^{\text{i}}]_{1,1}=\norm{(\mathbf{A+F}_1)^\top \mathbf{u}_1^{\text{i}}} $, the $k$-th iteration computes
\begin{equation}
\begin{aligned}
    &\bar{\mathbf{u}}_{k}^{\text{i}}=(\mathbf{A+E}_k)\mathbf{v}_k^{\text{i}}, & &\mathbf{u}^{\text{i}}=(\mathbf{I - U}_k^{\text{i}}(\mathbf{U}_k^{\text{i}})^\top)\bar{\mathbf{u}}_k^{\text{i}},   & &\mathbf{u}_{k+1}^{\text{i}}=\bfu^{\text{i}}/\norm{\bfu^{\text{i}}},\\
    &\bar{\mathbf{v}}_{k+1}^{\text{i}}=(\mathbf{A+F}_{k+1})^\top\mathbf{u}_{k+1}^{\text{i}},
    & &\mathbf{v}^{\text{i}} = (\mathbf{I - V}_k^{\text{i}}(\mathbf{V}_k^{\text{i}})^\top)\bar{\mathbf{v}}_{k+1}^{\text{i}},    & &\mathbf{v}_{k+1}^{\text{i}}=\mathbf{v}^{\text{i}}/\norm{\bfv^{\text{i}}},
\end{aligned}
\end{equation}
where $\mathbf{U}_k^{\text{i}} = [\mathbf{u}_1^{\text{i}},\dots,\mathbf{u}_k^{\text{i}}]\in \mathbb{R}^{m\times k}$ and $ \mathbf{V}_k^{\text{i}}=[\mathbf{v}_1^{\text{i}},\dots,\mathbf{v}_k^{\text{i}}] \in\mathbb{R}^{n\times k}$ are matrices with orthonormal columns.

After $k$ iterations, iGK algorithm computes upper Hessenberg matrix $\mathbf{M}_k^{\text{i}}\in\mathbb{R}^{(k+1)\times k} $ with $[\mathbf{M}_k^{\text{i}}]_{j,i}=(\mathbf{u}_j^{\text{i}})^\top\bar{\mathbf{u}}_i^{\text{i}} $ and $[\mathbf{M}_k^{\text{i}}]_{i+1,i}=\norm{\mathbf{u}^{\text{i}}} $ for $1\le j\le i\le k$ , lower triangular matrix $\mathbf{L}_{k+1}^{\text{i}}\in\mathbb{R}^{(k+1)\times(k+1) } $ with $[\mathbf{L}_{k+1}^{\text{i}}]_{i+1,j}=(\mathbf{v}_j^{\text{i}})^\top\bar{\mathbf{v}}_{i+1}^{\text{i}} $ and $[\mathbf{L}_{k+1}^{\text{i}}]_{i+1,i+1}=\norm{\mathbf{v}_i^{\text{i}}} $ for $1\le j\le i\le k$, with the following relationships 
\begin{equation}
    \begin{aligned}
 & (\mathbf{A}+\mathcal{E}_k)\mathbf{V}_k^{\text{i}} = \mathbf{U}_{k+1}^{\text{i}}\mathbf{M}_k^{\text{i}}, & \quad &\mathcal{E}_k=\sum_{i=1}^k \mathbf{E}_i\mathbf{v}_i^{\text{i}}(\mathbf{v}_i^{\text{i}})^\top, \\
&(\mathbf{A} + \mathcal{F}_{k+1})^\top \mathbf{U}_{k+1}^{\text{i}} = \mathbf{V}_{k+1}^{\text{i}}(\mathbf{L}_{k+1}^{\text{i}})^\top, & \quad &\mathcal{F}_{k+1} = \sum^{k+1}_{i=1} (\mathbf{u}_i^{\text{i}}(\mathbf{u}_i^{\text{i}})^\top)\mathbf{F}_i.
\end{aligned}
\end{equation}
To ensure the orthogonality of $\mathbf{V}_k^{\text{i}}$ and $\mathbf{U}_{k+1}^{\text{i}} $ under the inexactness of $\mathbf{A \text{ and }A^\top}$, iGK generates $\mathbf{M}_k^{\text{i}}$ and $\mathbf{L}_{k+1}^{\text{i}}$ instead of  $\mathbf{B}_k^{\text{s}} $ as in \cref{eq:GKB matrix fac}.

Here, the inexact LSQR (iLSQR) is used solve for \cref{eq:standard setup} with inexact $\bfA$ and $\bfA^\top$. At the $k$-th iteration, it computes
 \begin{equation}
\bfy_k^{\text{i}} = \argmin_{\bfy \in \R^k}  \norm{\bfM_k^{\text{i}} \bfy - \beta^{\text{i}} \bfe_1}_2^2
\end{equation}
and obtain the solution $\bfx_k^{\text{i}} = \bfV_k^{\text{i}} \bfy_k^{\text{i}}$. Note that due to inxactness in the MVPs, the $k$-th iteration of iLSQR does not minimize the exact residual $\norm{\bfA\bfx_k - \bfb}$ along $\bfx_k^{\text{i}} \in \text{span}\{ \bfV_k^{\text{i}}\}$, and span\{$\mathbf{V}_k ^{\text{i}}$\} and span\{$\mathbf{U}_{k+1}^{\text{i}}$\} are no longer Krylov subspaces.

\subsection{Generalized Golub-Kahan Bidiagonalization}

Generalized Golub-Kahan (genGK) Bidiagonalization is an iterative method designed to solve linear systems in the generalized least-squares sense \cite{Arioli2013generalized, ChungGeneralized}. Referring to the problem we have in \cref{eq:reformulated MAP} where exact MVPs with $\bfA$ and $\bfA^\top$ could be achieved, genGk generates two sets of orthogonal vectors that span the Krylov subspaces 
$\mathcal{K}_k(\mathbf{A^\top R}^{-1} \mathbf{AQ, A^\top R}^{-1}\mathbf{b})$ and $\mathcal{K}_k(\mathbf{AQA^\top R}^{-1}, \mathbf{b}) $. 

Given matrices $\mathbf{A,\ R,\ Q}$ and vector $\mathbf{b}$, the genGK process initiates $\mathbf{u}_1^{\text{g}} =\mathbf{b}/\beta_1^{\text{g}} $ where 
$\beta_1^{\text{g}} = \norm{\mathbf{b}}_{\mathbf{R}^{-1}}$,  and $\alpha_1^{\text{g}}\mathbf{v}_1^{\text{g}}=\mathbf{A^\top R}^{-1} \mathbf{u}_1^{\text{g}}$ where $\alpha_1^{\text{g}} = \norm{ \mathbf{v}_1^{\text{g}}}_\mathbf{Q}$. We use the superscript `g' to denote the vectors and matrices computed from the genGK process. At the $k$-th iteration, genGK generates $\mathbf{u}_{k+1}^{\text{g}}$ and $\mathbf{v}_{k+1}^{\text{g}}$ through
\begin{equation}
\begin{aligned}
     &\beta_{k+1}^{\text{g}}\mathbf{u}_{k+1}^{\text{g}} = \mathbf{AQv}_{k+1}^{\text{g}} -\alpha_k^{\text{g}}\mathbf{u}_k^{\text{g}},  \\
     &\alpha_{k+1}^{\text{g}}\mathbf{v}_{k+1}^{\text{g}} = \mathbf{A^\top R}^{-1}\mathbf{u}_{k+1}^{\text{g}} - \alpha_{k+1}^{\text{g}}\mathbf{v}_k^{\text{g}},
\end{aligned}
\end{equation}
where scalars $\alpha_k^{\text{g}}, \beta_k^{\text{g}} \ge 0$ are chosen such that $\norm{\mathbf{u}_k^{\text{g}}}_{\mathbf{R}^{-1}}= \norm{\mathbf{v}_k^{\text{g}}}_{\mathbf{Q}} =1 $. 

After $k$ iterations, the algorithm generates
\[
\mathbf{B}_k^{\text{g}} = \begin{bmatrix}
    \alpha_1^{\text{g}} &0&\dots&0\\
    \beta_2^{\text{g}} & \alpha_2^{\text{g}}&\ddots&\vdots  \\
    0& \beta_3^{\text{g}} &\ddots&0\\
    \vdots&\ddots &\ddots & \alpha_k^{\text{g}}\\
    0 &\dots &0 & \beta_{k+1}^{\text{g}}
 \end{bmatrix}, \quad \mathbf{U}_{k+1}^{\text{g}} = [\mathbf{u}^{\text{g}},\dots,\mathbf{u}_{k+1}^{\text{g}} ],\quad \mathbf{V}_k^{\text{g}}=[\mathbf{v}_1^{\text{g}},\dots,\mathbf{v}_k^{\text{g}}]
\] with the following relations holding up to machine precision
\begin{equation}
\begin{aligned}
   \mathbf{U}_{k+1}^{\text{g}}\beta_1^{\text{g}}\mathbf{e}_1 &= \mathbf{b}, \\
    \mathbf{AQV}_k^{\text{g}} &= \mathbf{U}_{k+1}^{\text{g}}\mathbf{B}_k^{\text{g}},\\
   \mathbf{ A^\top R}^{-1} \mathbf{U}_{k+1}^{\text{g}} &= \mathbf{V}_k^{\text{g}}(\mathbf{B}_k^{\text{g}})^\top +\alpha_{k+1}^{\text{g}}\mathbf{v}_{k+1}^{\text{g}}\mathbf{e}_{k+1}^\top ,
\end{aligned}
\end{equation}
and matrices $\mathbf{U}_{k+1}^{\text{g}} $ and $\mathbf{V}_k^{\text{g}} $ satisfies the following orthogonality conditions:
\begin{equation}
    (\mathbf{V}_k^{\text{g}})^\top \mathbf{QV}_k^{\text{g}}=\mathbf{I}_k  \quad\text{and} \quad (\mathbf{U}_{k+1}^{\text{g}})^\top \mathbf{R}^{-1}\mathbf{U}_{k+1}^{\text{g}}=\mathbf{I}_{k+1}.
\end{equation}

Again, we can consider an iterative method where we seek a solution $\bfx_k^{\text{g}}$ in the span of $\bfV_k^{\text{g}}$, i.e., $\bfx_k^{\text{g}} = \bfV_k^{\text{g}} \bfy_k^{\text{g}} $, where at each iteration, we solve for $\bfy_k^{\text{g}}$
\begin{equation}
\bfy_k^{\text{g}} = \min_{\bfy \in \R^n} \norm{\bfB_k^{\text{g}}\bfy - \beta_1^{\text{g}}\bfe_1 }_2^2 +\frac{\lambda^2}{2} \norm{\bfy}_2^2.
\end{equation}
We call this the genLSQR method.

\section{Iterative methods based on inexact generalized Golub-Kahan decomposition}
\label{sec:Inexact genGK}

This section proposes a iterative solver based on the inexact generalized Golub-Kahan (igenGK) decomposition. This method aims to solve problems where matrices $\mathbf{A}\text{ and } \mathbf{Q}$ are really large that can only be accessed through MVPs via function evaluations, for example the least squares problem as stated in \cref{eq:reformulated MAP}.  The main difference is that we allow potential inexactness in $\bfA$ and $\bfA\t$.

\subsection{Inexact generalized Golub-Kahan decomposition} \label{sec:inexact A}

Assume the MVPs with $\mathbf{A}$ and $\mathbf{A}^\top$ cannot be performed exactly, following the setup as in \cref{eq:inexact A def}, and the covariance kernel matrix $\bfQ$ can only be accessed through MVPs. Here, we propose an inexact generalized Golub-Kahan decomposition method for solving such problem.

The igenGK decomposition method combines the  iGK and genGK approaches. With initializations  $\mathbf{u}_1=\mathbf{b}/\beta$ where $\beta = \norm{\mathbf{b}}_{\mathbf{R}^{-1}} $, $\mathbf{v}_1=(\mathbf{A+F}_1)^\top\mathbf{R}^{-1} \mathbf{u}_1 /[\mathbf{L}]_{1,1} $ where $[\mathbf{L}]_{1,1}=\norm{(\mathbf{A+F}_1)^\top\mathbf{R}^{-1} \mathbf{u}_1}_\mathbf{Q} $, at the $k$-th iteration, it computes
\begin{equation*}
    \begin{aligned}
    &\bar{\mathbf{u}}_{k}=(\mathbf{A+E}_k)\mathbf{Q} \mathbf{v}_k, & &\mathbf{u}=(\mathbf{I - U}_k\mathbf{U}_k^\top\mathbf{R}^{-1} )\bar{\mathbf{u}}_k,   & &\mathbf{u}_{k+1}=\mathbf{u/\norm{u}}_{\mathbf{R}^{-1}} ,\\
    &\bar{\mathbf{v}}_{k+1}=(\mathbf{A+F}_{k+1})^\top \mathbf{R}^{-1} \mathbf{u}_{k+1},
    & &\mathbf{v} = (\mathbf{I - V}_k\mathbf{V}_k^\top\mathbf{Q} )\bar{\mathbf{v}}_{k+1},    & &\mathbf{v}_{k+1}=\mathbf{v/\norm{v}}_\mathbf{Q},
    \end{aligned}
\end{equation*}
where $\mathbf{U}_k = [\mathbf{u}_1,\dots,\mathbf{u}_k]\in \mathbb{R}^{m\times k}$ and $ \mathbf{V}_k=[\mathbf{v}_1,\dots,\mathbf{v}_k] \in\mathbb{R}^{n\times k}$ are matrices with orthonormal columns, separately with respect to $\bfR^{-1}$ and $\bfQ$:
\begin{equation}
\mathbf{U}_{k+1}^\top \mathbf{R}^{-1}\mathbf{U}_{k+1}=\mathbf{I}_{k+1} \quad\text{and}  \quad \mathbf{V}_k^\top \mathbf{QV}_k=\mathbf{I}_k .
\end{equation}

Thus, after $k$ iterations, it generates an upper Hessenberg matrix $\mathbf{M}_k\in\mathbb{R}^{(k+1)\times k}$ with $[\mathbf{M}_k]_{j,i}=\mathbf{u}_j^\top \mathbf{R}^{-1} \bar{\mathbf{u}}_i $ and $[\mathbf{M}_k]_{i+1,i}=\norm{\mathbf{u}}_{\mathbf{R}^{-1} } $ for $1\le j\le i\le k$; a lower triangular matrix $\mathbf{L}_{k+1}\in\mathbb{R}^{(k+1)\times(k+1) } $  with $[\mathbf{L}_{k+1}]_{i+1,j}=\mathbf{v}_j^\top\mathbf{Q} \bar{\mathbf{v}}_{i+1} $ and $[\mathbf{L}_{k+1}]_{i+1,i+1}=\norm{\mathbf{v}}_\mathbf{Q} $ for $1\le j\le i\le k$; as well as following relationships
\begin{equation}
    \begin{aligned}
    (\mathbf{A}+\mathcal{E}_k)\mathbf{Q}\mathbf{V}_k &= \mathbf{U}_{k+1}\mathbf{M}_k, \\
        (\mathbf{A} +\mathcal{F}_{k+1})^\top\mathbf{R}^{-1}\mathbf{U}_{k+1} &= \mathbf{V}_{k+1}\mathbf{L}_{k+1}^\top,
    \end{aligned}
\end{equation}
where $\mathcal{E}_k=\sum_{i=1}^k \mathbf{E}_i\mathbf{Q}\mathbf{v}_i\mathbf{v}_i^\top$ , $\mathcal{F}_{k} = \sum_{i=1}^{k} (\mathbf{u}_i\mathbf{u}_i^\top) \mathbf{R}^{-1}\mathbf{F}_i$. To see this, notice that
\begin{equation*}
    \begin{aligned}
 \relax[(\mathbf{A}+\mathbf{E}_1)\mathbf{Q}\mathbf{v}_1 ,\ldots ,(\mathbf{A}+\mathbf{E}_k)\mathbf{Q}\mathbf{v}_k] =&\mathbf{A}\mathbf{QV}_k+[\mathbf{E}_1\mathbf{Q}\mathbf{v}_1, \dots,  \mathbf{E}_k\mathbf{Q}\mathbf{v}_k]\mathbf{V}_k^\top \mathbf{QV}_k\ \quad \\
    =&\mathbf{A}\mathbf{Q}\mathbf{V}_k + \left(\sum_{i=1}^k \mathbf{E}_i\mathbf{Q}\mathbf{v}_i\mathbf{v}_i^\top \right)\mathbf{Q}\mathbf{V}_k\\
    =&(\mathbf{A}+\mathcal{E}_k)\mathbf{Q}\mathbf{V}_k,
\end{aligned}
\end{equation*}
and similarly,
\begin{equation*}
\begin{aligned}
     \relax[(\mathbf{A}+\mathbf{F}_1) ^\top \mathbf{R}^{-1}\mathbf{u}_1,  \dots ,(\mathbf{A}+\mathbf{F}_{k})^\top \mathbf{R}^{-1}\mathbf{u}_{k}]
    =&(\mathbf{A} +\mathcal{F}_{k})^\top\mathbf{R}^{-1}\mathbf{U}_{k}.
\end{aligned}
\end{equation*}

The following algorithm summarizes the inexact generalized Golub-Kahan decomposition process. Given the inexactness of forward matrix $\mathbf{A}$, the algorithm takes input $\mathbf{A}$ as a function operator and performs $\widehat{\bfA}_k \bfx$ and $\widetilde{\bfA}^\top_k \bfy$ for the $k$-th iteration MVP, $\mathbf{V}_{k-1}\in\mathbb{R}^{n \times (k-1)}$ with $\mathbf{V}_0$ as an empty matrix, $\mathbf{U}_k\in\mathbb{R}^{m\times k}$ with $\mathbf{U}_1=\mathbf{b}/\norm{\mathbf{b}}_{\mathbf{R}^{-1}} $, $\mathbf{M}_{k-1}\in\mathbb{R}^{k\times (k-1)}$ an upper Hessenberg matrix, and $\mathbf{C}_{k} = \mathbf{L}_k^\top \in\mathbb{R}^{k\times k}$ an upper triangular matrix.

\begin{algorithm}[H]
\caption{igenGK decomposition with Reorthogonalization}
\label{alg:igenGk}
\begin{algorithmic}[1]
\STATE{\textbf{Input:} $\mathbf{ A,\mathbf{Q},R},\mathbf{V}_{k-1},\mathbf{U}_k,\mathbf{C}_{k-1}, \mathbf{M}_{k}$ }
\IF {$k=1$}
    {\STATE{ $c_{kk}\mathbf{v}_k = \widetilde{\bfA}_k^\top\mathbf{R}^{-1}\mathbf{u}_k\text{ where } c_{kk}=\norm{\mathbf{v}_k}_\mathbf{Q}$}
    }
\ELSE{ 
    \STATE $\mathbf{v} =\widetilde{\bfA}_k^\top\mathbf{R}^{-1}\mathbf{u}_k$ }

    \FOR{$j=1,\dots,k-1$}
    \STATE {$\mathbf{v} =\mathbf{v}-c_{ik}\mathbf{v}_j \text{ where } c_{ik} = \mathbf{v}_j^\top\mathbf{Q}\mathbf{v}$}
    \ENDFOR
    
    \STATE{$\mathbf{v}_k=\mathbf{v}/c_{kk} \text{ where } c_{kk} = \norm{\mathbf{v} }_\mathbf{Q}$}
    \STATE{$\mathbf{u} =  \widehat{\bfA}_k\mathbf{Qv}_k$}
    \FOR{$j=1,\dots,k$}
    \STATE{ $\mathbf{u = u}-m_{jk}\mathbf{u_j} \text{ where } m_{jk}=\mathbf{u}_j^\top\mathbf{R}^{-1}\mathbf{u} $}
\ENDFOR
\STATE{ $\mathbf{u}_{k+1} = \mathbf{u}/m_{k+1,k} 
\text{ where }m_{k+1,k}=\mathbf{\norm{u}_{R^{-1}} }$
}
\ENDIF
\RETURN $\mathbf{V}_k=[\mathbf{V}_{k-1}\ \mathbf{v}_k], \mathbf{U}_{k+1}=[\mathbf{U}_k\ \mathbf{u}_{k+1}],  \mathbf{M}_{k+1}, \mathbf{C}_{k} $
\end{algorithmic}
\end{algorithm}

\subsection{Solving the LS problem}
In the subsection above, we introduce the igenGK process as an iterative method to generate a subspace for the solution. Further, we would like to solve the least-squares problem  \cref{eq:reformulated MAP} through a sequence of projected LS problems which we denote the inexact generalized LSQR  (igenLSQR) method. In particular, at the $k$-th step, we seek solution $\mathbf{x}_k\in\text{span}\{\mathbf{V}_k\}$, and define residual $\bfr_k \equiv (\bfA + \mathcal{E}_k)\bfQ\bfx_k - \bfb$, then we obtain the following igenLSQR problem,
\begin{equation}
\min_{\mathbf{x}_k\in \mathrm{span}\{\mathbf{V}_k\}} 
    \frac{1}{2} \norm{(\mathbf{A+\mathcal{E}}_k)\mathbf{Qx}_k-\mathbf{b}}_{\mathbf{R}^{-1}}^2 + \frac{\lambda^2}{2} \norm{\mathbf{x}_k}^2_\mathbf{Q}.
    \label{eq:proj LS}
\end{equation}
For now, we may assume $\lambda$ is fixed. Given $\mathbf{x}_k=\mathbf{V}_k\bfy_k$ where $\bfy_k \in\mathbb{R}^k$, we have
\[
(\mathbf{A}+\mathcal{E}_k)\mathbf{QV}_k\bfy_k-\mathbf{b} = \mathbf{U}_{k+1}\mathbf{M}_k\bfy_k - \mathbf{b}=\mathbf{U}_{k+1}(\mathbf{M}_k\bfy_k -\beta \mathbf{e}_1)
\]
and
$$\bfy_k ^\top \mathbf{V}_k^\top \mathbf{QV}_k\bfy_k = \bfy_k^\top \bfy_k.$$ 
The above LS problem \cref{eq:proj LS} could be formulated as follows
\begin{equation}
    \min_{{\bfy_k}\in\mathbb{R}^k} \frac{1}{2}\norm{(\mathbf{A}+\mathcal{E}_k)\mathbf{QV}_k\bfy_k-\mathbf{b}}^2_{\bfR^{-1}}+\frac{\lambda^2}{2} \norm{\mathbf{V}_k\bfy_k}^2_\mathbf{Q} 
\end{equation}
which is equivalent to
\begin{equation}
\min_{{\bfy_k}\in\mathbb{R}^k} \frac{1}{2}\norm{\mathbf{M}_k\bfy_k-\beta \mathbf{e}_1}^2_2+\frac{\lambda^2}{2} \norm{\bfy_k}^2_2 .
    \label{eq: proj LS_z}
\end{equation}
This approach is stemming from the standard LSQR  \cite{paige1982lsqr, paige1982lsqr2} and generalized LSQR \cite{ChungGeneralized}.

After obtaining a solution for the projected problem, the solution $\bfs_k$ for the original problem \cref{eq:MAP} could be recovered as
\begin{equation}
    \bfs_k = \boldsymbol{\mu} + \bfQ\bfx_k = \boldsymbol{\mu} + \bfQ\bfV_k\bfy_k.
\end{equation}

\subsection{Inexact generalized hybrid approach}
\label{sec:hybrid}
So far, we have assumed the Tikhonov regularization parameter $\lambda$ is known a priori. However, in practice, obtaining a good regularization parameter is crucial but could also be difficult, especially when the problem is large in scale. If a poor $\lambda$ is chosen, it may lead to an imbalance between the residual and perturbation error, thus ending up with poor solutions. In this section, we propose an inexact generalized hybrid approach, where the regularization parameter 
can be automatically estimated at each iteration. 

This method follows from many previous works, where the problem is first projected down to a lower dimensional space and the projected problem \cref{eq: proj LS_z} would be further solved through various parameter selection methods \cite{hanke1993regularization, hansen2010discrete, chung2024computational}.

While various regularization parameter methods are available \cite{chung2008weighted, chung2015hybrid, kilmer2001choosing, renaut2010regularization}, here we consider the discrepancy principle (DP) and the weighted generalized cross validation (WGCV) approach. To provide a benchmark for parameter selection method comparisons, both will be compared with the optimal approach, where the regularization parameter $\lambda_{\text{opt}} $ is chosen to minimize the 2-norm of the error between the reconstructed solution and the true solution
\begin{equation}
\lambda_{\text{opt}} = \argmin_\lambda \norm{\bfs_k(\lambda) - \bfs_{\text{true}}}_2^2,
\label{eq:opt parameter}
\end{equation} 
where $\bfs_\text{true}$ denotes the true solution and $\bfs_k(\lambda)$ denotes the solution computed at the $k$-th iteration using regularization parameter $\lambda$.

At the $k$-th iteration, we define the residual  as
$$\bfr_k (\lambda) = \bfA \bfs_k(\lambda) - \bfd = \bfM_k\bfy_k(\lambda) - \beta_1\bfe_1. $$

The discrepancy principle (DP) method chooses $\lambda = \lambda_{\text{DP} }$ to minimize the distance between the residual norm and the level of noise in the observation
\begin{equation}
\lambda_\text{DP} =  \argmin_\lambda \left| \norm{\bfM_k\bfy_k(\lambda) - \beta_1\bfe_1}_2 - \nu_\text{DP} \norm{\boldsymbol{\epsilon}}_2 \right|
\end{equation}  
where $\nu_\text{DP}$ is a user chosen constant. 

The weighted generalized cross validation (WGCV) method is another common approach for selecting regularization parameters when the level of noise is unknown. This method follows from the statistical technique cross validation. By arbitrarily leaving out one element of the observed data $\bfd$, cross validation aims to find a good regularization parameter which is able to  predict the missing element. However, these approaches can be expensive, so variants such as the generalized cross validation \cite{golub1965calculating} method have been considered. Developed in the context of hybrid projection methods, the WGCV method selects the regularization parameter as 

\begin{equation}
    \lambda_{\text{WGCV}} = \argmin_\lambda \frac{\norm{\bfM_k\bfy_k(\lambda) - \beta\bfe_1 }_2^2}{\text{trace}\left(\bfI_k - \omega \bfM_k(\bfM_k^\top \bfM_k + \lambda^2 \bfI)^{-1} \bfM_k\right)^2},
\end{equation}
where $\omega$ is a weight parameter that can be user defined or estimated during the iterative process \cite{chung2008weighted,chung2024computational}.

\section{Numerical Experiments}
\label{sec:numerical experiment}

\begin{figure}[tbhp]
    \centering
    \includegraphics[width=0.9\linewidth]{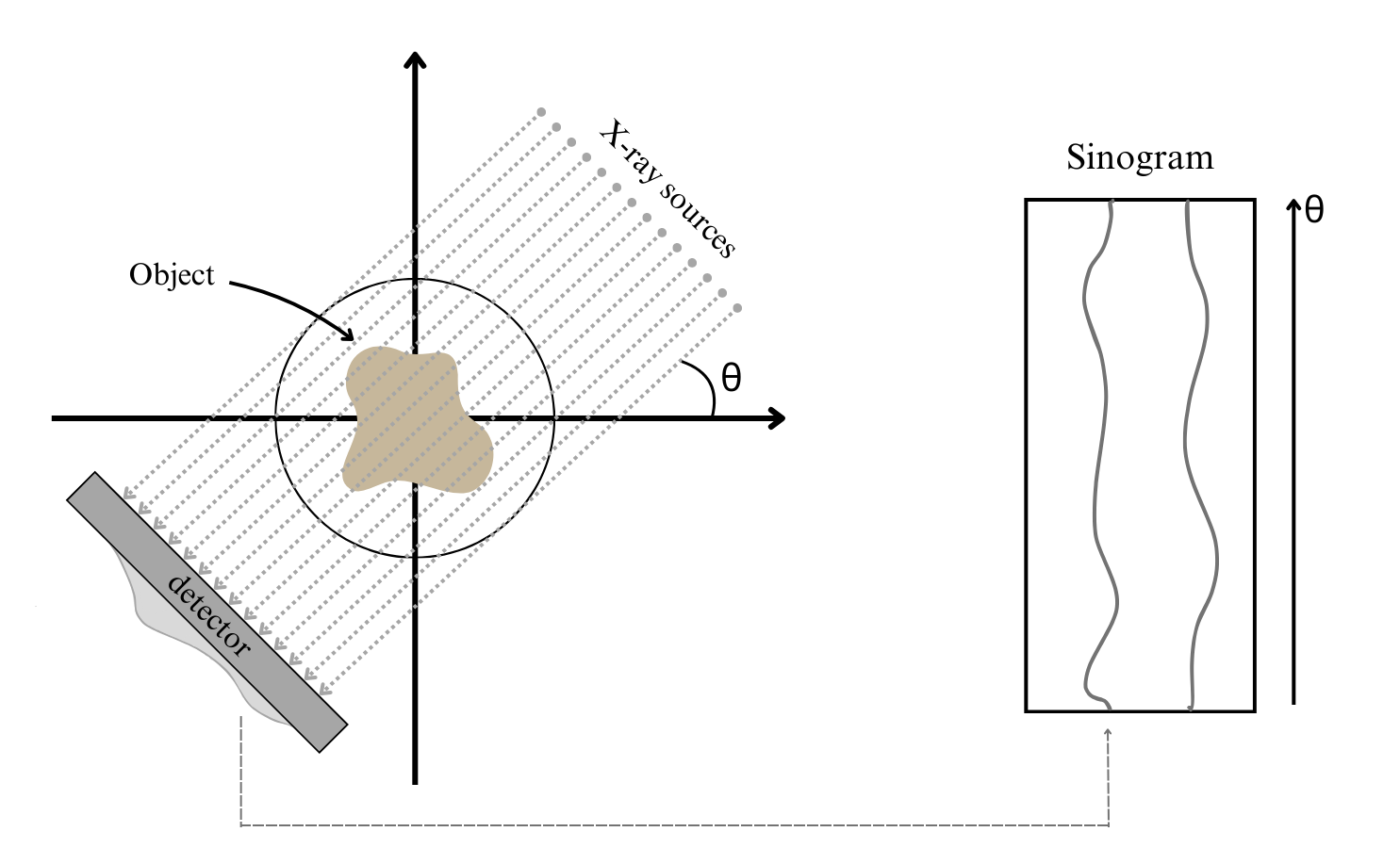}
    \caption{Set up of parallel-beam X-ray CT where the detector collects data from a finite number of projection angles $\boldsymbol{\theta}$ and the data from a set of $\boldsymbol{\theta}$ constitute the sinogram.}
    \label{fig: tomo}
\end{figure}

In this section, we consider a X-ray computed tomographic (CT) reconstruction problem, where the goal is to reconstruct the cross-section of an object from data collected along the X-rays penetrating the object. As shown in \cref{fig: tomo}, the X-ray source-detector pair rotates $180^\circ$ around the object and a finite number of projections are collected at degrees $\boldsymbol{ \theta}$, where $\boldsymbol{\theta}$ are projection angles and the resulting observation constitutes a 2D sinogram image which could be further vectorized as observation $\bfd$ in the inverse problem.

All experiments below are performed in MATLAB R2023a, using the `PRtomo' example from the IRTools package \cite{gazzola2019ir}. Within each experiment, the proposed igenGK decomposition is compared with GK, iGK, and genGK with respect to their performance in image reconstruction. Also the relative reconstruction error norm at each iteration $k$ is computed as
$$e_k = \norm{\bfs_k - \bfs_\text{true}}_2 / \norm{\bfs_\text{true}}_2.$$

\subsection{Numerical experiment 1: Comparison of iterative methods without regularization\label{sec:num exp no regularization A}}

We begin with an experiment that compares the  performance of GK, iGK, genGK, and igenGK methods without using any additional regularization on the projected problem (i.e., $\lambda = 0$).

PRtomo generates a `medical' phantom image of size $128\times 128$, so we stack the columns together to obtain $\bfs_\text{true} \in\mathbb{R}^{16384}$. The goal of this problem is to reconstruct the true image through the given forward X-ray CT process $\bfA\in \mathbb{R}^{6516\times16384}$  and observation $\bfd \in \mathbb{R}^{6516}$. Suppose the matrix $\bfA$ is modeled after a parallel-beam process where the projections are recorded at every $5^\circ$ with $\boldsymbol{\theta}=[1,6,11,\dots, 176]$.  Further, we assume the true image $\bfs \sim \mathcal{N}(\mathbf{0}, \lambda^{-2}\bfQ )$ where Matérn covariance kernel $\bfQ$ is defined by $\nu = 1.5,\ \ell=0.01$, and Gaussian white noise is added to $\bfd_\text{true}$ to make it more  realistic, so $\bfd = \bfd_\text{true}+\boldsymbol{\epsilon}$ where $\boldsymbol{\epsilon} \sim \mathcal{N}(\mathbf{0}, \sigma\bfI) $, $\sigma$ is chosen such that $\norm{\boldsymbol{\epsilon}}_2 / \norm{\bfd_\text{true}}_2 = 0.04$ indicating $4\%$ noise level.

To simulate errors in the forward model, we incorporate random inexactness within each MVP with $\bfA$ and $\bfA^\top$. To implement this, the forward matrix is defined as an object (i.e., using object-oriented programming) and is accessed through function evaluations for MVPs with $\bfA$ and $\bfA^\top$. More precisely, at each iteration $k$, $\widehat{\bfA}_k\bfx = (\bfA + \bfE_k)\bfx$ and $\widetilde{\bfA}_k^\top\bfx = (\bfA^\top + \bfF_k)\bfx$ and $\bfE_k , \bfF_k$ are matrices whose entries are random numbers generated i.i.d. from the Gaussian distribution $\mathcal{N}(\bf0, \beta)$ with $\beta>0$.

Firstly, we would like to verify the following relationships for the igenGK decomposition:
\begin{align}
    \mathbf{A}^\top \mathbf{U}_k&= \mathbf{V}_{k}\mathbf{L}_{k}^\top, \label{relation:1} \\
    \mathbf{A}\mathbf{Q}\mathbf{V}_k& = \mathbf{U}_{k+1} \bfM_k, \label{relation:2}\\
    \mathbf{I}_k = \mathbf{V}_{k}^\top \mathbf{Q}\mathbf{V}_{k} \quad &\text{and}\quad \mathbf{I}_{k+1} = \mathbf{U}_{k+1}^\top \mathbf{R}^{-1} \mathbf{U}_{k+1}.\label{relation:3}
\end{align}
In \cref{tab:igenGK A relations}, we present the results for various degrees of inexactness.  We consider $k=50$ and the following values for $\beta$: $10^{-2}, 10^{-4},$ and $10^{-6}$, and we provide the four relations using the metric 
\[
\text{Error} = \frac{\text{norm of left hand side $-$ right hand side}}{\text{norm of left hand side}},
\]
where the Frobenius norm is used.
Results indicate that the discrepancies between the left and right hand sides of \cref{relation:1} and \cref{relation:2} are proportional to the amount of random inexactness introduced to each MVP with $\mathbf{A}$ and $\mathbf{A}^\top$, while \cref{relation:3} holds up to machine precision. 

\begin{table}[tbhp]
\centering
\begin{tabular}{@{}llll@{}}
\toprule
\multicolumn{1}{c}{\multirow{2}{*}{igenGK Relations}} & \multicolumn{3}{c}{Error}                                                            \\ \cmidrule(l){2-4} 

\multicolumn{1}{c}{}    & $\beta=10^{-2}$ & $\beta=10^{-4}$ & $\beta=10^{-6}$ \\ \midrule
$\norm{\mathbf{A}^\top \mathbf{U}_k-\mathbf{V}_{k}\mathbf{L}_{k}^\top}/\norm{\mathbf{A}^\top \mathbf{U}_k}$
& 5.26e-02                   & 5.26e-04                   & 5.26e-06                   \\
$\norm{\mathbf{A}\mathbf{Q}\mathbf{V}_k - \mathbf{U}_{k+1} \bfM_k} / \norm{\mathbf{A}\mathbf{Q}\mathbf{V}_k}$ & 3.05e-02                 & 3.07e-04                   & 3.07e-06                   \\
$ \norm{\mathbf{V}_{k}^\top \mathbf{Q}\mathbf{V}_{k} - \mathbf{I}_k} / \sqrt{k}$ & 1.86e-15                   & 2.63e-15                   & 1.43e-15                   \\
$\norm{\mathbf{U}_{k+1}^\top \mathbf{R}^{-1} \mathbf{U}_{k+1} - \mathbf{I}_{k+1}} /  \sqrt{k+1}$ & 1.64e-14                   & 1.03e-15                   & 1.06e-14                   \\ \bottomrule
\end{tabular}
\caption{Verifications of igenGK relationships for varying degrees of inexactness.}
\label{tab:igenGK A relations}
\end{table}

To further evaluate its capability in constructing the solution subspace, we 
evaluate the image reconstruction performance by directly solving $\bfy_k$ from \cref{eq: proj LS_z} with $\lambda=0$, i.e. the projected problems are solved without additional regularization technique. Then, we derive the MAP estimate $\mathbf{s}_k$ in \cref{eq:MAP} by computing $\mathbf{s}_k=\mathbf{Qx}_k=\mathbf{Q}\mathbf{V}_{k}\bfy_k$.

\begin{figure}[tbhp]
\centering
\begin{subfigure}{.5\textwidth}
  \centering
  \includegraphics[width=\linewidth]{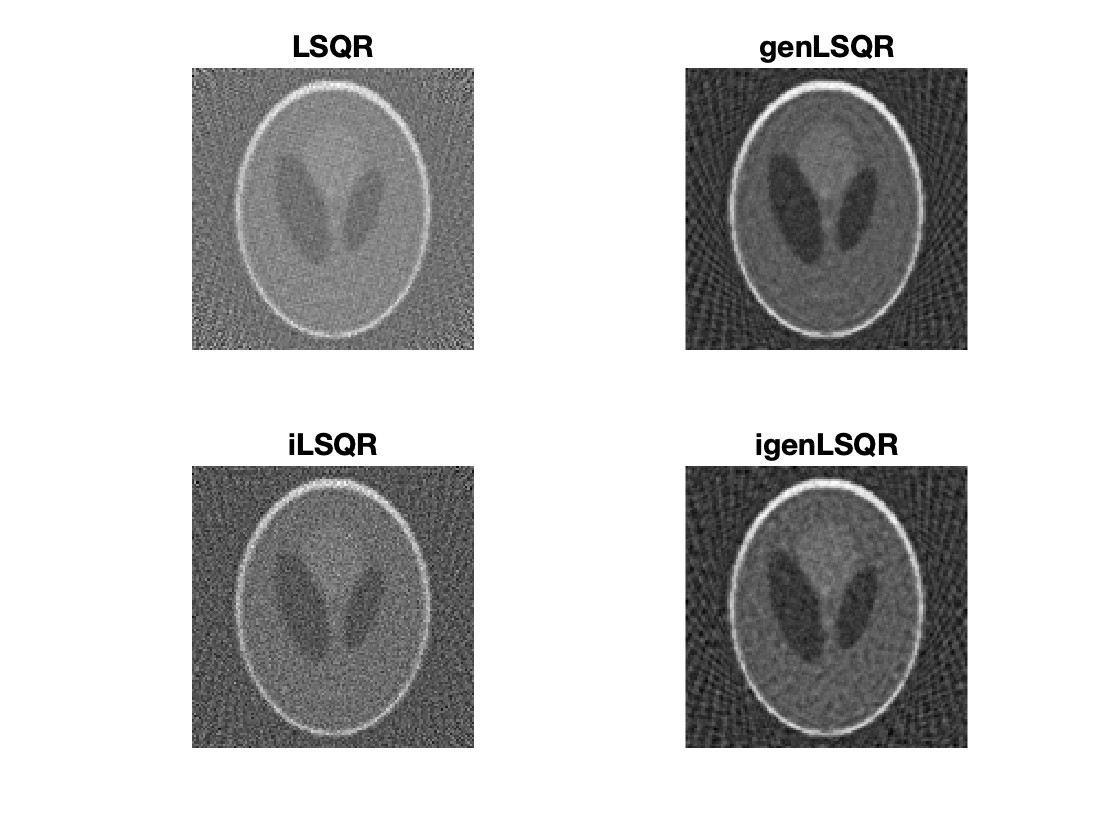}
  \caption{Reconstruction}
  \label{fig:one-step comp A}
\end{subfigure}%
\begin{subfigure}{.5\textwidth}
  \centering
  \includegraphics[width=\linewidth]{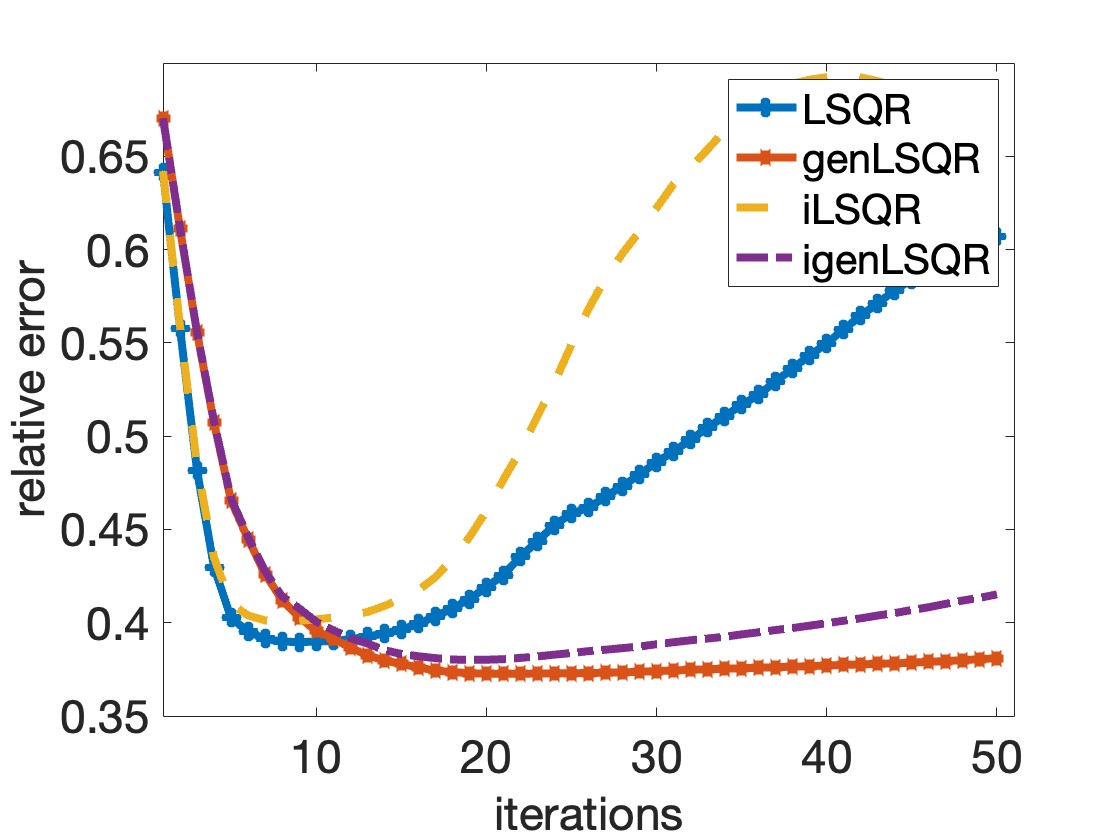}
  \caption{Relative reconstruction error norms}
  \label{fig:one-step e-norm comp A}
\end{subfigure}
\caption{Both (a) and (b) present comparisons of standard,  generalized, inexact, and inexact generalized iterative methods for image reconstruction. The legend LSQR corresponds to the standard method, genLSQR corresponds to the generalized method, iLSQR cooresponds to the inexact method, and igenLSQR correponds to the inexact generalized method.  In (a) we present the image reconstructions at the end of $50$ iterations, and in (b) we provide the relative errors norms along the iterations.}

\end{figure}

For comparative analysis, we further compute the MAP estimates using the GK, iGK, and genGK iterative methods. For genGK and igenGK, which are generalized methods that require selecting a prior covariance matrix $\bfQ$, we set $\bfQ$ to be a kernel matrix defined by the Matérn kernel with $\nu=1.5$ and $\ell = 0.01$. Both iGK and igenGK are inexact methods, and we use the inexact forward matrix $\bfA$ with $\beta = 10^{-2}$. The image reconstructions obtained from each method are presented in Figure \ref{fig:one-step comp A}, while Figure \ref{fig:one-step e-norm comp A} presents the relative error norms along iterations. The legend ``LSQR" corresponds to standard GK iterations,``genLSQR" corresponds to the genGK iterative method,``iLSQR" corresponds to the inexact iterative method iGK, and ``igenGK" correponds to inexact generalized iterative method igenGK.

We observe that, without additional regularization, all methods exhibit semiconvergence: they firstly converge in the first few iterations but after which the
error increases dramatically when the inverted noise starts to dominate the solution \cite{hanke1993regularization}. So the reconstructed images at the end of $50$ iterations are not desirable since they are contaminated with noise. From both \cref{fig:one-step comp A} and \cref{fig:one-step e-norm comp A}, while the solution from LSQR is close to iLSQR and genLSQR is close to igenLSQR, the generalized methods start with slower convergence but eventually perform better at later iterations with lower relative errors and less noisy images. Further, the methods that involve inexactness (iLSQR and igenLSQR) perform worse but not significantly worse than their counterparts (LSQR and genLSQR), which is promising.

\subsection{Numerical experiment 2: Comparisons of hybrid methods with optimal regularization}

As we observed in \cref{sec:num exp no regularization A}, this large-scale problem without additional regularization on projected problems demonstrates semiconvergence. Thus, we seek the parameter-choice methods which automatically select the regularization parameters within the projected solution space at each iteration  to achieve more stable and reliable solutions. In this experiment, following the same setup as the previous one,  we further present the performance of the igenGK method when adopting a hybrid approach. Here, we consider the optimal regularization parameter $\lambda = \lambda_{\text{opt}} $  from \cref{eq:opt parameter}.

\begin{figure}[tbhp]
\centering
\begin{subfigure}{.5\textwidth}
  \centering
  \includegraphics[width=\linewidth]{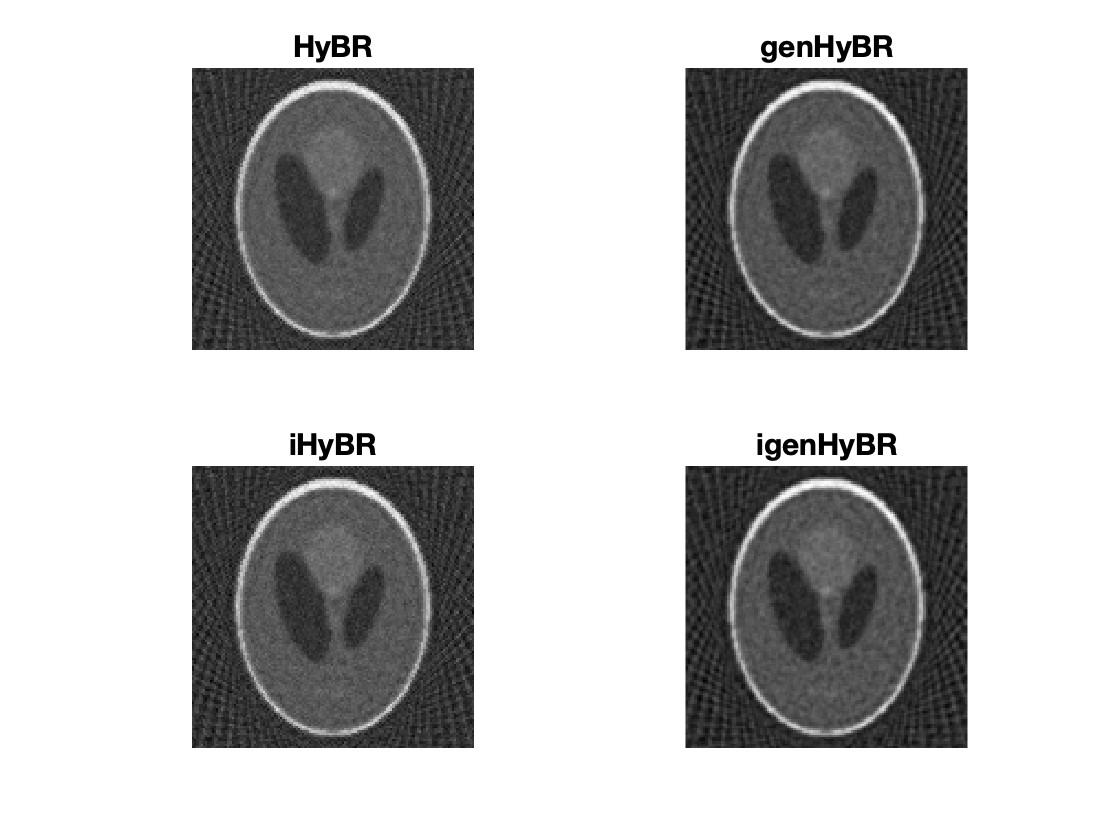}
  \caption{Reconstruction}
  \label{fig: hybrid reconstruct A}
\end{subfigure}%
\begin{subfigure}{.5\textwidth}
  \centering
  \includegraphics[width=\linewidth]{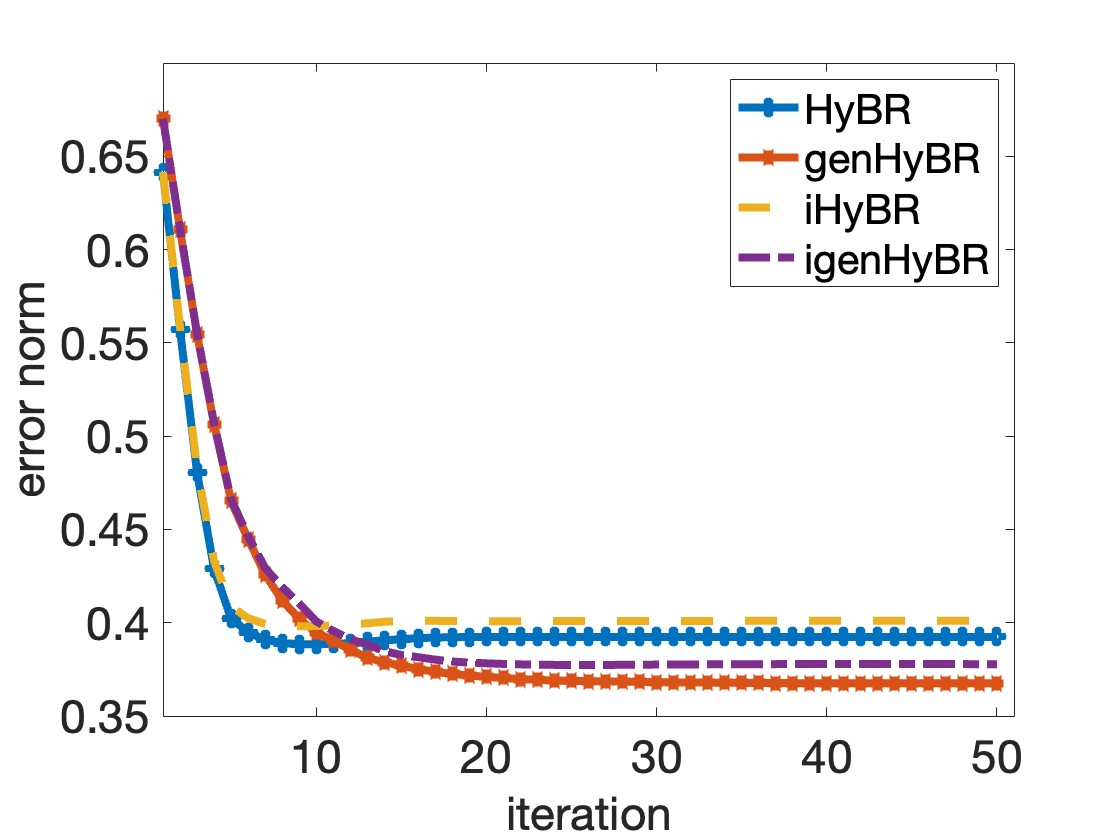}
  \caption{Relative reconstruction error norms}
  \label{fig: hybrid error A}
\end{subfigure}
\caption{Comparisons of standard, generalized, inexact, and inexact generalized hybrid approach with optimal regularization parameter. HyBR corresponds to the hybrid method based on GK, genHyBR corresponds to the generalized hybrid method based on genGK, iHyBR corresponds to the inexact hybrid method based on iGK, and igenHyBR corresponds to the inexact generalized hybrid method based on igenGK. In (a) we present the image reconstructions at the end of $50$ iterations, (b) provides the relative errors norms along the iterations.}

\end{figure}

Figure \ref{fig: hybrid reconstruct A} and Figure \ref{fig: hybrid error A} separately show the image reconstruction results and relative errors for all four iterative methods in adopting the hybrid approach.``HyBR" corresponds to standard hybrid method based on GK, ``genHyBR" corresponds to the generalized hybrid method, ``iHyBR" corresponds to the inexact hybrid method, and ``igenHyBR" corresponds to the inexact generalized hybrid method.  All of these results correspond to using the optimal regularization parameter.

We observe that the hybrid method not only leads to better reconstructed images but also stabilizes the convergence. Again, the igenHyBR we propose achieves a low relative reconstruction error, being only slightly higher than its exact counterpart.

\subsection{Numerical experiment 3: Comparison of regularization parameters}

While the optimal regularization parameter is the parameter corresponding to the Tikhonov reconstruction that is closest to the true solution in the 2-norm sense, it is not realistic in  applications since it requires the true solution which is not available. Therefore, we further compare two additional methods for choosing the regularization parameter, as introduced in \cref{sec:hybrid}, that are more practical and commonly used.

\begin{figure}[tbhp]
\centering
\includegraphics[width=10cm]{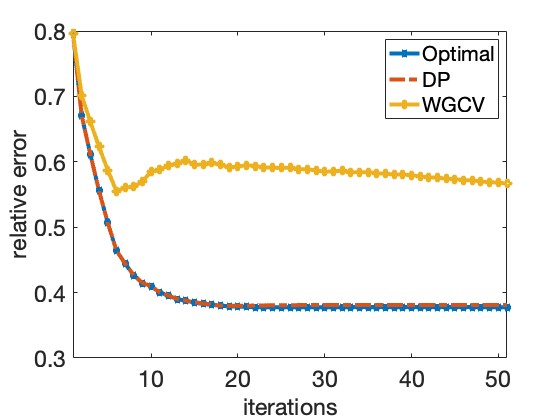}
\caption{Comparison of optimal, DP, and WGCV methods in choosing the Tikhonov regularization parameter $\lambda$. The relative error norms are computed along iterations where each regularization method is incorporated into igenHyBR in computing solutions.}
\label{fig: compare parameter}
\end{figure}

In the following, we present an experiment to compare the different regularization parameter selection methods concerning their efficacy in approximating the true solution. We adopt the igenHyBR method with regularization parameters selected optimally, using the discrepancy principle (DP), and using the weighted generalized cross validation regularization (WGCV) method, and we make comparisons of these three hybrid approaches.  We see in \cref{fig: compare parameter} that results for DP are similar to those for the optimal regularization parameter. This is very promising since the DP method does not require knowledge of true solution; however, it does require a noise level estimate. For these results, we use the actual noise level. WGCV seems to struggle with selecting a good regularization parameter, but we observed similar results for the genHyBR method and additional tuning of the weight parameter could improve the solutions.

\subsection{Numerical experiment 4: Inexactness in projection angles}

In this experiment, we consider a more realistic scenario where
the inexactness of the forward process arises from inaccurate projection angles. Continuing with the X-ray CT example, the forward process is determined by numerous factors, such as the ray type (whether it uses parallel beam, fan beam, or cone beam), the projections angles (the degree angles chosen to conduct projections), and the number and spread of beams emitted from the X-ray sources at each projection angle. Uncertainty in any of these setups could lead to large variances in the  observed sinogram and thus make the reconstructed image undesirable. 

In reality, the degree at which the X-rays are taken and recorded is one of the main sources of uncertainty in CT reconstruction. Though we assume the exact projection angles are set before the scanning process, error arises due to inaccurate calibration of the equipment or unforeseen incidents, such as the involuntary movement of the scanned object, or the external factors causing vibrations in X-ray sources and detector \cite{uribe2022hybrid}.

The following example illustrates a scenario where the inexactness in the projection angles decreases over iterations. This is reasonable if one considers integrating these inexact methods within a larger optimization scheme where the angles are being updated and improved.  That is, each solution could be adopted at the next iteration to better improve the estimated angles at which the X-rays are collected. Given the the accurate projection angles $\boldsymbol{\theta}_{\text{true}} = [1,6,11,\dots,176]^\top$, we simulate this scenario by using at each iteration $k$ a set of inexact projection angles $\boldsymbol{\theta}_k = \boldsymbol{\theta}_{\text{true}}+ \alpha_k \bfe_k $ where $\alpha_k$ decreases logarithmically along iterations and $\bfe_k$ is a random vector drawn from a standard normal distribution for each $k$.

\begin{figure}[tbhp]
 \begin{subfigure}{0.49\textwidth}
     \centering
     \includegraphics[width=\textwidth]{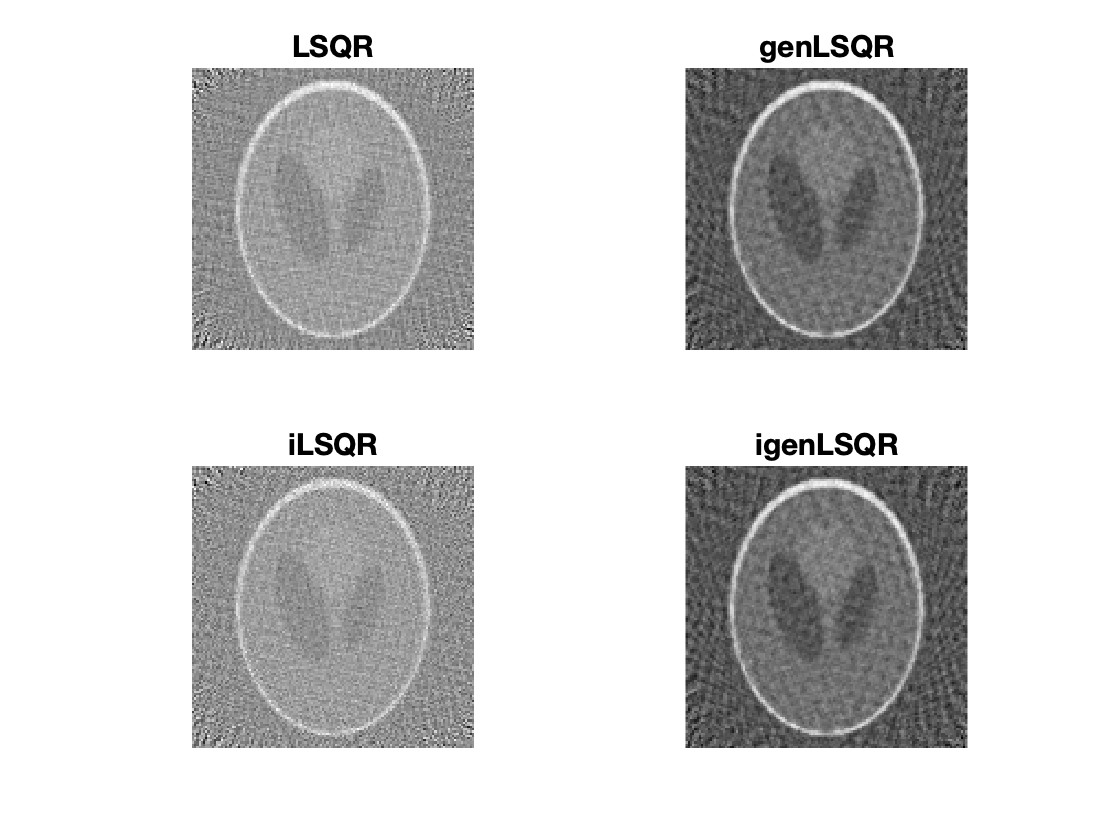}
 \end{subfigure}
 \hfill
 \begin{subfigure}{0.49\textwidth}
     \centering
     \includegraphics[width=\textwidth]{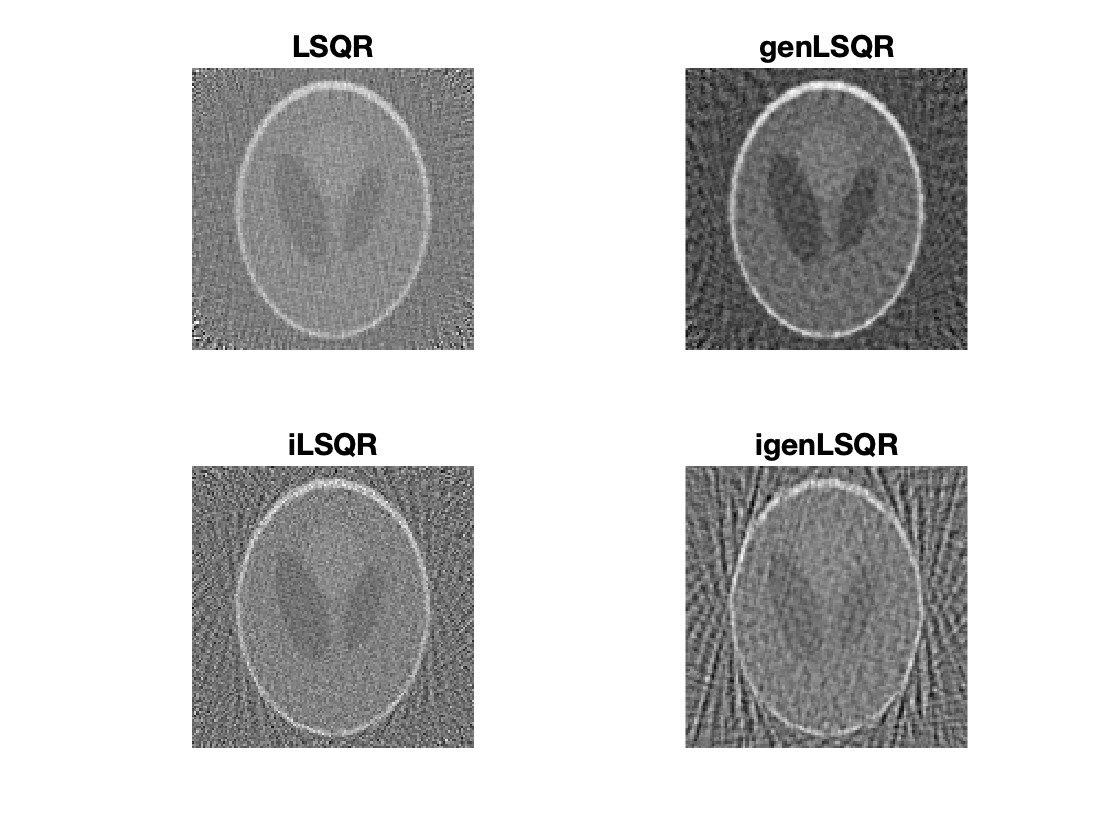}
 \end{subfigure}
 
 \medskip
 \begin{subfigure}{0.49\textwidth}
 \centering
     \includegraphics[width=\textwidth]{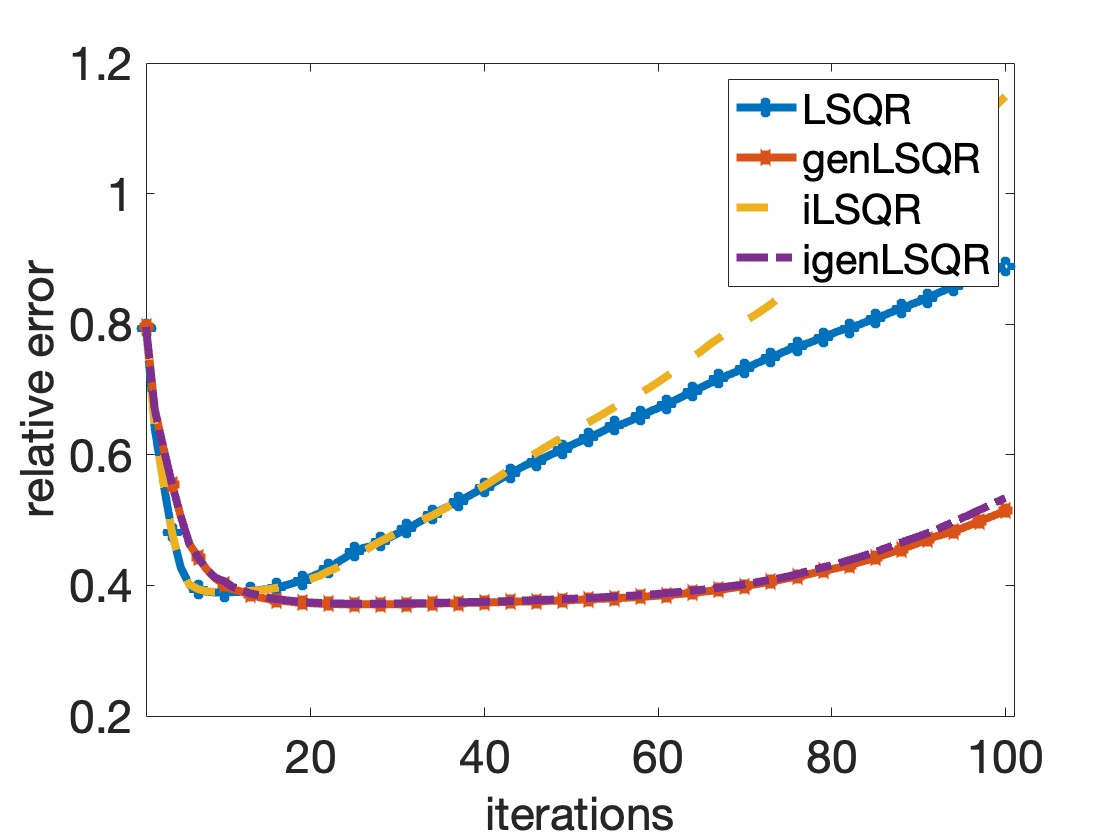}
     \caption{Small initial inexactness}
     \label{fig:small inx angle}
 \end{subfigure}
 \hfill
 \begin{subfigure}{0.49\textwidth}
     \centering
     \includegraphics[width=\textwidth]{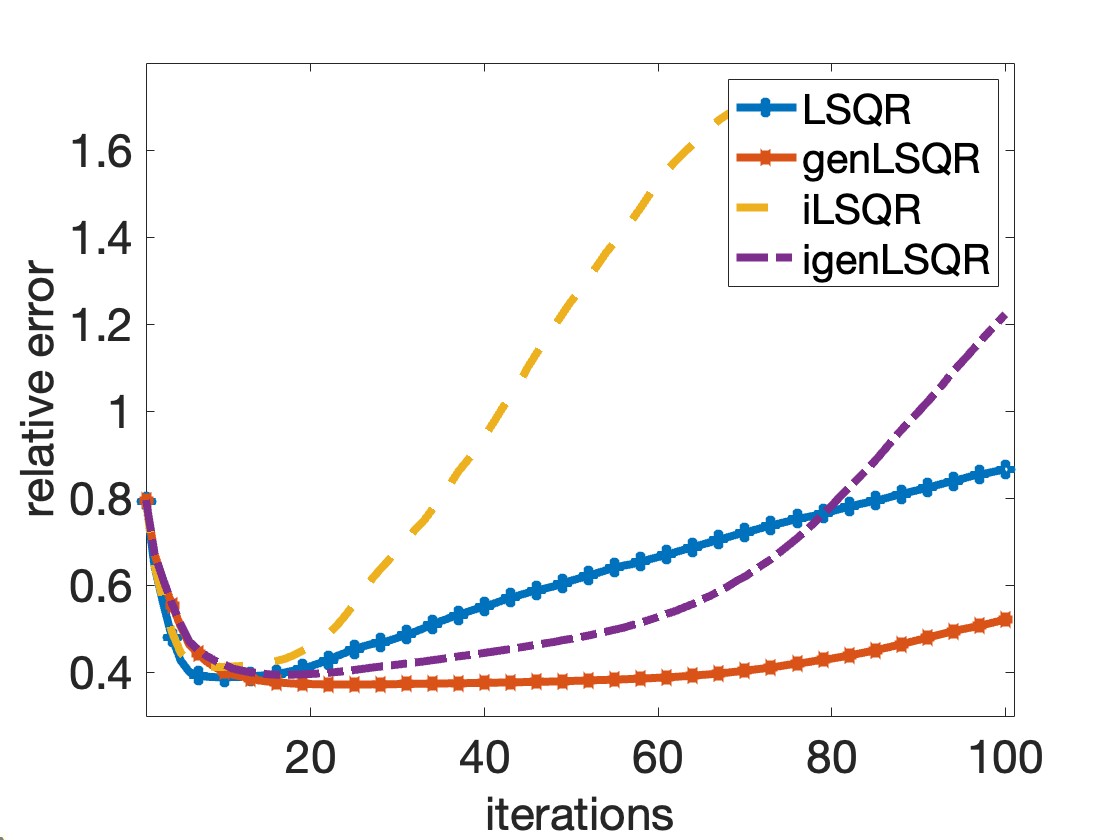}
     \caption{Large initial inexactness}
     \label{fig:large inx angle}
 \end{subfigure}

 \caption{Comparisons of reconstruction results for sinograms acquired at inexact angles, solved without regularization. (a) starts from smaller degree of inexactness ($\alpha_1 = 10^{-1}$), (b) starts from larger degree of inexactness ($\alpha_1 = 1$).}
 \label{fig:inx angle no regu}

\end{figure}

In \cref{fig:inx angle no regu}, we present the image reconstruction and relative error results for two cases: the first case is presented in \cref{fig:small inx angle}  where we start with a small degree of inexactness $\alpha_1= 10^{-1}$ and decreases logarithmically down to $\alpha_{100} = 10^{-6}$,  and in the second case as shown in \cref{fig:large inx angle}, we start from a larger inexactness $\alpha_1 = 10^{0}$ and similarly decreases down to $10^{-6}$. Both of them are solved without additional regularization on the projected problems ($\lambda=0$) and we extend the iteration number to 100 to fully show the semiconvergence in relative error curves.

\begin{figure}[tbhp]
 \begin{subfigure}{0.49\textwidth}
     \centering
     \includegraphics[width=\textwidth]{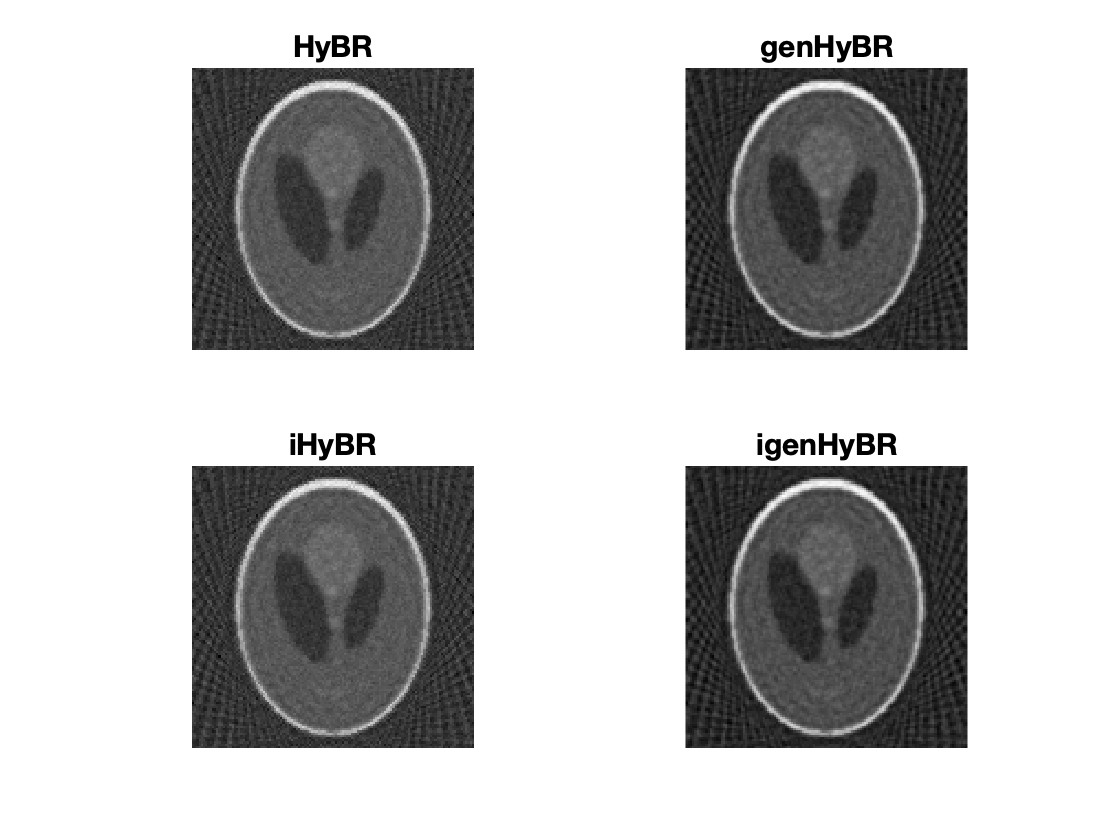}
 \end{subfigure}
 \hfill
 \begin{subfigure}{0.49\textwidth}
     \centering
     \includegraphics[width=\textwidth]{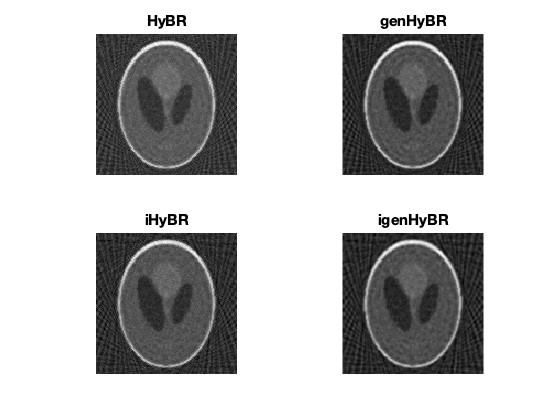}
 \end{subfigure}
 
 \medskip
 \begin{subfigure}{0.49\textwidth}
 \centering
     \includegraphics[width=\textwidth]{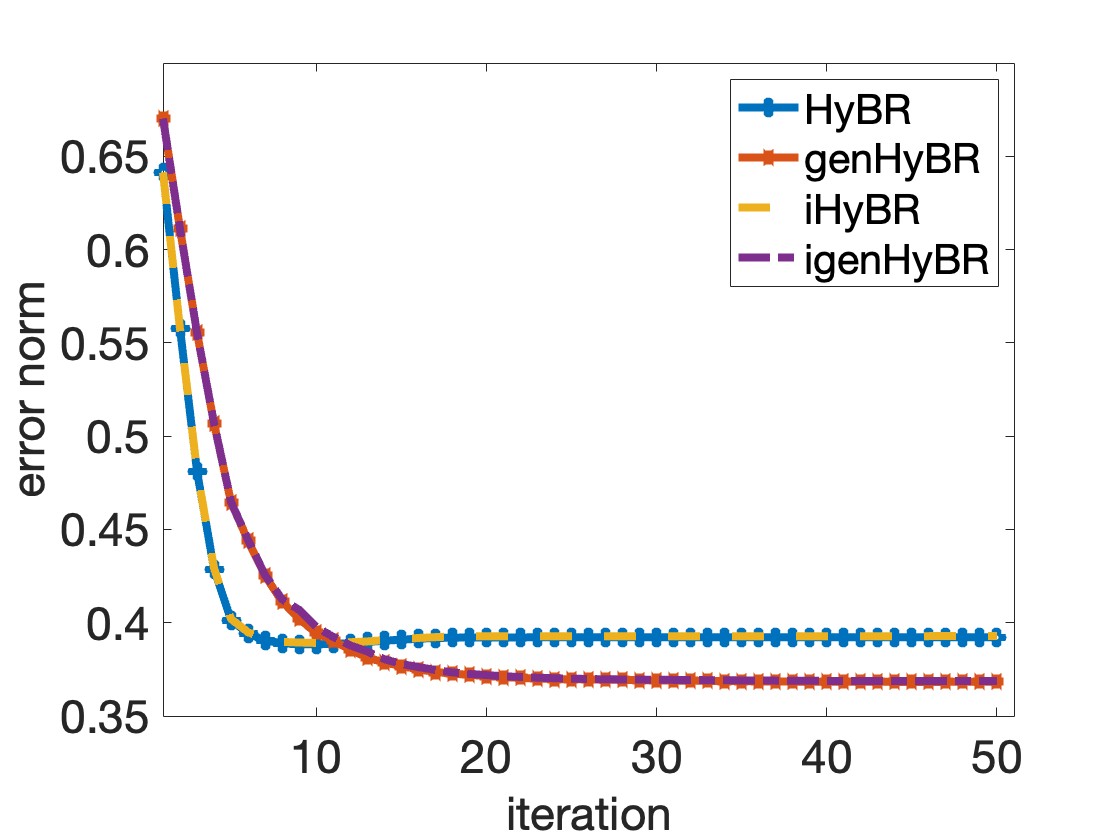}
     \caption{Small initial inexactness}
     \label{fig:small inx angle hybrid}
 \end{subfigure}
 \hfill
 \begin{subfigure}{0.49\textwidth}
     \centering
     \includegraphics[width=\textwidth]{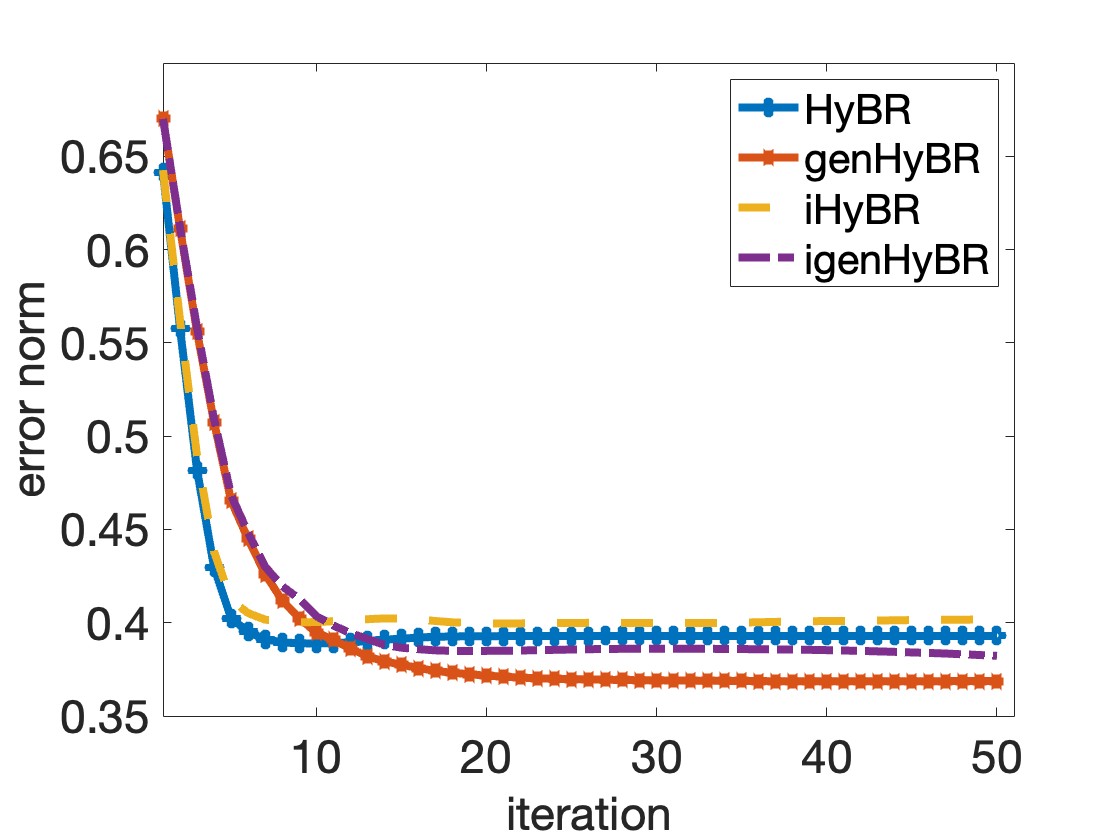}
     \caption{Large initial inexactness}
     \label{fig:large inx angle hybrid}
 \end{subfigure}

 \caption{Comparisons of reconstruction results for sinograms acquired at inexact angles, solved through hybrid approach with optimal regularization. (a) starts from smaller degree of inexactness ($\alpha_1 = 10^{-1}$), (b) starts from larger degree of inexactness ($\alpha_1 = 1$).}
 \label{fig:inx angle opt regu}

\end{figure}

Further incorporating regularization on the projected problems with optimal regularization parameters, we present the same two examples solved using igenHyBR in \cref{fig:inx angle opt regu}. Here, we observe better reconstructed images and more stable relative error curves. However, though both cases end at the same degree of inexactness with $\alpha_{50} = 10^{-6}$, igenHyBR in \cref{fig:small inx angle hybrid} achieves better performance than in \cref{fig:large inx angle hybrid}, with the relative error curve overlapping with its exact counterpart. This indicates that the starting degree of inexactness is decisive in reconstruction performance.

\section{Conclusion}
\label{sec:conclusion}
In this paper, we developed an inexact generalized hybrid iterative method for efficiently computing solutions to large-scale Bayesian inverse problems with inexactness in the forward process. Unlike previous approaches which assume the exact MVPs with the forward matrix are achievable, our method adapts to the inexactness and allows effective and efficient computations. Moreover, we developed a hybrid iterative projection method that combines the igenGK projection approach with Tikhonov regularization on the projected problem. Compared to approaches that use the exact forward model, the inexact genGK methods demonstrate stability and accuracy, even for problems with large degrees of inexactness.
Numerical results on tomographic image reconstruction problems validate the effectiveness of the approach and show its adaptability to real-world conditions, such as inaccurate projection angles in CT imaging.

\section*{Acknowledgments}
The author received the Wilson Family Undergraduate STEM Research Award and would like to thank the donors.
Also, this work was partially supported by the National Science Foundation program under grant DMS-2411197. Any opinions, findings, and conclusions or recommendations expressed in this material are those of the author(s) and do not necessarily reflect the views of the National Science Foundation.
\appendix
\section{Choice of covariance kernel}
\label{sec:cov kernel}
We choose the covariance kernel to be from Matérn family for two main reasons. First, it allows for varying levels of smoothness and  correlation structures between points. 
Given the entries of a covariance matrix are computed as $\mathbf{Q}_{ij} = \kappa(\mathbf{x}_i,\mathbf{x}_j) $, where $\{\mathbf{x}\}_{i=1}^n $ are spatial points in domain. We define the  covariance matrix coming from Matérn family to be
\begin{equation}  \kappa(\mathbf{x}_i,\mathbf{x}_j)=C_{\alpha, \nu}(r) = \frac{1}{2^{\nu-1}\Gamma(\nu)} \left(\sqrt{2\nu} \alpha r \right)^\nu K_\nu\left(\sqrt{2\nu} \alpha r\right)
\end{equation}
where $r=\norm{\mathbf{x}_i,\mathbf{x}_j}_2 $, $\Gamma$ is the Gamma function, $K_\nu(\cdot)$ is the modified Bessel function of the second kind of order $\nu$, and $\alpha$ is a scaling factor. 
When $\nu = 1/2$, $C_{\alpha, \nu}$ corresponds to the exponential covariance function, and when $\nu \rightarrow \infty$, $C_{\alpha, \nu}$ converges to the Gaussian covariance function.

Second, this choice of $\bfQ$ allows for efficient storage and MVPs. Normally, Matérn covariance matrix is very dense and thus expensive to store and compute. Using naive approach would costs $\mathcal{O}(n^2)$ for both storing and performing one MVP. However, as one of transnational (or stationary) invariant covariance kernels, the cost per MVP with Matérn family of covariance kernels could be reduced to $\mathcal{O}(n\log n)$ by exploiting methods such as the fast Fourier transform (FFT) or the fast multipole method \cite{nowak2003efficient}. This project uses the connection between FFT and Toeplitz (1D) /block-Toeplitz (2D) structure for efficient computation of $\mathbf{Qx}$.

\bibliographystyle{siamplain}
\bibliography{siuro_article}

\end{document}